\documentclass[a4paper,12pt]{amsart}
\usepackage{amsfonts}
\usepackage{amssymb}
\usepackage{ifthen}
\usepackage{graphicx}
\usepackage{color}
\nonstopmode \numberwithin{equation}{section}
\newtheorem{thm}{Theorem}
\newtheorem{lem}{Lemma}
\newtheorem{cor}{Corollary}[section]
\newtheorem{prop}{Proposition}[section]

\newtheorem{cl}{Claim}
\newtheorem{ca}{Case}
\newtheorem{sca}{Subcase}
\newtheorem{scl}{Subclaim}
\newtheorem{conj}{Conjecture}

\theoremstyle{definition}
\newtheorem{defn}{Definition}

\newtheorem{op}[equation]{Open Problem}
\newtheorem{ques}[equation]{Question}
\newtheorem{rem}{Remark}[section]
\newtheorem{exam}[equation]{Example}

\newcounter {own}
\def\theown {\thesection       .\arabic{own}}

\newenvironment{pf}[1][]{%
 \vskip 3mm
 \noindent
 \ifthenelse{\equal{#1}{}}%
  {{\slshape Proof. }}%
  {{\slshape #1.} }%
 }%
{\qed\bigskip}

\newcounter{alphabet}
\newcounter{tmp}
\newenvironment{Thm}[1][]{\refstepcounter{alphabet}%
\bigskip%
\noindent%
{\bf Theorem \Alph{alphabet}}%
\ifthenelse{\equal{#1}{}}{}{ (#1)}%
{\bf .} \itshape}{\vskip 8pt}

\makeatletter
\newcommand{\Ref}[1]{\@ifundefined{r@#1}{}{\setcounter{tmp}{\ref{#1}}\Alph{tmp}}}
\makeatother

\newenvironment{Lem}[1][]{\refstepcounter{alphabet}%
\bigskip%
\noindent%
{\bf Lemma \Alph{alphabet}}%
{\bf .} \itshape}{\vskip 8pt}

\newcommand{\ID}{{\mathbb D}}
\newcommand{\IB}{{\mathbb B}}

\def\be{\begin{equation}}
\def\ee{\end{equation}}

\newcommand{\bee}{\begin{enumerate}}
\newcommand{\eee}{\end{enumerate}}

\newcommand{\blem}{\begin{lem}}
\newcommand{\elem}{\end{lem}}
\newcommand{\bthm}{\begin{thm}}
\newcommand{\ethm}{\end{thm}}
\newcommand{\bcor}{\begin{cor}}
\newcommand{\ecor}{\end{cor}}
\newcommand{\beg}{\begin{exam}}
\newcommand{\eeg}{\end{exam}}
\newcommand{\begs}{\begin{examples}}
\newcommand{\eegs}{\end{examples}}
\newcommand{\bdefe}{\begin{defn}}
\newcommand{\edefe}{\end{defn}}
\newcommand{\bprob}{\begin{prob}}
\newcommand{\eprob}{\end{prob}}
\newcommand{\bques}{\begin{ques}}
\newcommand{\eques}{\end{ques}}
\newcommand{\bei}{\begin{itemize}}
\newcommand{\eei}{\end{itemize}}
\newcommand{\bcon}{\begin{conj}}
\newcommand{\econ}{\end{conj}}
\newcommand{\bop}{\begin{op}}
\newcommand{\eop}{\end{op}}

\newcommand{\bca}{\begin{ca}}
\newcommand{\eca}{\end{ca}}
\newcommand{\bsca}{\begin{sca}}
\newcommand{\esca}{\end{sca}}

\newcommand{\bcl}{\begin{cl}}
\newcommand{\ecl}{\end{cl}}

\newcommand{\bscl}{\begin{scl}}
\newcommand{\escl}{\end{scl}}

\newcommand{\bcons}{\begin{conjs}}
\newcommand{\econs}{\end{conjs}}
\newcommand{\bprop}{\begin{propo}}
\newcommand{\eprop}{\end{propo}}
\newcommand{\br}{\begin{rem}}
\newcommand{\er}{\end{rem}}
\newcommand{\brs}{\begin{rems}}
\newcommand{\ers}{\end{rems}}
\newcommand{\bo}{\begin{obser}}
\newcommand{\eo}{\end{obser}}
\newcommand{\bos}{\begin{obsers}}
\newcommand{\eos}{\end{obsers}}
\newcommand{\bpf}{\begin{pf}}
\newcommand{\epf}{\end{pf}}
\newcommand{\ba}{\begin{array}}
\newcommand{\ea}{\end{array}}
\newcommand{\beq}{\begin{eqnarray}}
\newcommand{\beqq}{\begin{eqnarray*}}
\newcommand{\eeq}{\end{eqnarray}}
\newcommand{\eeqq}{\end{eqnarray*}}

\newcommand{\ra}{\rightarrow}

\newcommand{\ds}{\displaystyle}

\newcounter{minutes}
\divide\time by 60
\newcounter{hours}
\multiply\time by 60 \addtocounter{minutes}{-\time}

\begin{document}

\bibliographystyle{amsplain}
\title{Stable geometric properties of pluriharmonic and biholomorphic mappings, and Landau-Bloch's theorem
}

\def\thefootnote{}
\footnotetext{ \texttt{\tiny File:~\jobname .tex,
          printed: \number\day-\number\month-\number\year,
          \thehours.\ifnum\theminutes<10{0}\fi\theminutes}
} \makeatletter\def\thefootnote{\@arabic\c@footnote}\makeatother

\author{Sh. Chen }
\address{Sh. Chen, Department of Mathematics and Computational
Science, Hengyang Normal University, Hengyang, Hunan 421008,
People's Republic of China.} \email{mathechen@126.com}

\author{S. Ponnusamy $^\dagger $
}
\address{S. Ponnusamy,
Indian Statistical Institute (ISI), Chennai Centre, SETS (Society
for Electronic Transactions and security), MGR Knowledge City, CIT
Campus, Taramani, Chennai 600 113, India. }
\email{samy@isichennai.res.in, samy@iitm.ac.in}

\author{X. Wang 
}

\address{X. Wang, Department of Mathematics,
Hunan Normal University, Changsha, Hunan 410081, People's Republic
of China.} \email{xtwang@hunnu.edu.cn}

\keywords{Pluriharmonic   mapping, linearly connected domain,  Landau-Bloch's theorem.\\
$^\dagger$ {\tt This author is on leave from the
Department of Mathematics,
Indian Institute of Technology Madras, Chennai-600 036, India.\\
}
}


\begin{abstract}
In this paper, we investigate some properties of pluriharmonic
mappings defined in the unit ball. First, we discuss some geometric
univalence criteria on pluriharmonic mappings, and then establish a
Landau-Bloch theorem for a class of pluriharmonic mappings.
\end{abstract}


\maketitle \pagestyle{myheadings} \markboth{ Sh. Chen, S. Ponnusamy
and X. Wang }{Stable pluriharmonic and biholomorphic mappings}

\section{Introduction and main results  }\label{csw-sec1}
As usual, $\mathbb{C}^{n}$ denotes the complex Euclidean space of
$n$ variables $z=(z_{1},\ldots,z_{n})\in \mathbb{C}^{n}$ with the
standard {\it Hermitian inner product}  $\langle z,w\rangle :=
\sum_{k=1}^nz_k\overline{w}_k$ and norm $ \|z\|:={\langle
z,z\rangle}^{1/2},$ where $w=(w_{1},\ldots,w_{n})$, and
$\overline{w}_k$ $(1\leq k \leq n)$ denotes the complex conjugate of
$w_k$ with
$\overline{w}=(\overline{w}_{1},\ldots,\overline{w}_{n})$. For $a\in
\mathbb{C}^n$,
$$\IB^n(a,r)=\{z\in \mathbb{C}^{n}:\, \|z-a\|<r\}
$$
denotes the (open) ball of radius $r>0$ with center $a$ and
$$\partial \IB^n(a,r)=\{z\in \mathbb{C}^{n}:\, \|z-a\|=r\}.
$$
Also, we let $\IB^n(r):=\IB^n(0,r)$, and use $\IB^n$ to denote the unit ball $\IB^n(1)$, and
$\mathbb{D}=\mathbb{B}^1$.


Throughout, $\mathcal{H}(\IB^n, \mathbb{C}^n)$ denotes the set of
all  mappings $f$ from $\mathbb{B}^{n}$ into $\mathbb{C}^{n}$ which
are continuously differentiable as mappings into $\mathbb{R}^{2n}$
with $f=(f_{1},\ldots,f_{n})$ and $f_{j}(z)=u_{j}(z)+iv_{j}(z)$
($1\leq j\leq n$), where $u_{j}$ and $v_{j}$  are real-valued
functions from $\mathbb{B}^{n}$ into $\mathbb{R}$.

For a complex-valued and differentiable function $f$ from
$\mathbb{B}^{n}$ into $\mathbb{C}$, we introduce (see for instance
\cite{CPW-1,CPW6,CPW-2})
$$\nabla f=\left(\frac{\partial f}{\partial z_{1}},\ldots,\frac{\partial f}{\partial z_{n}}
\right )\;\;\mbox{and}\;\; \overline{\nabla}f =\left (\frac{\partial
f}{\partial \overline{z}_{1}},\ldots,\frac{\partial f}{\partial
\overline{z}_{n}}\right ).
$$

  For $f\in \mathcal{H}(\IB^n, \mathbb{C}^n)$, we use  $J_{f}$ to
  denote
the {\it real Jacobian matrix} of $f$ (cf. \cite{HG1}). 
Moreover, for each $f=(f_{1},\ldots,f_{n})\in \mathcal{H}(\IB^n,
\mathbb{C}^n)$, denote by
$$Df=\big(\nabla f_1,\ldots,\nabla f_n\big)^T
$$
the matrix whose rows are the complex gradients $\nabla
f_1,\ldots,\nabla f_n$, and let
$$\overline{D}f=\big(\overline{\nabla}f_1,\ldots,\overline{\nabla}f_n\big)^T,
$$
where $T$ means the matrix transpose.


For an $n\times n$ complex matrix $A$, we introduce the {\it
operator norm}
$$\|A\|=\sup_{z\neq0}\frac{\|Az\|}{\|z\|}=\max\left\{\|A\theta\|:\, \theta\in\partial\mathbb{B}^{n}\right\}.
$$
We use $L(\mathbb{C}^{n},\mathbb{C}^{m})$ to denote the space of
continuous {\it linear operators} from $\mathbb{C}^{n}$ into
$\mathbb{C}^{m}$ with the operator norm, and let $I_{n}$ be the {\it
identity operator} in $L(\mathbb{C}^{n},\mathbb{C}^{n})$.


It follows from \cite[Theorem 4.4.9]{R} that a real-valued function $u$
defined on a simply connected domain $G$ is pluriharmonic if and only if $u$ is the real part
of a holomorphic function on $G$. Clearly, a mapping
$f:\,\mathbb{B}^{n}\ra \mathbb{C}$ is pluriharmonic if and only if
$f$ has a representation $f=h+\overline{g}$, where $g$ and $h$ are
holomorphic. We refer to \cite{HG1,CPW-1,DHK,I,R} for the definition and further
details on pluriharmonic mappings.

A {\it vector-valued mapping} $f \in\mathcal{H}(\IB^n,
\mathbb{C}^n)$ is said to be  pluriharmonic, if each of its
component functions is a  pluriharmonic mapping from $\mathbb{B}^{n}$
into $\mathbb{C}$. We denote by $\mathcal{PH}(\mathbb{B}^{n},
\mathbb{C}^n)$ the set of all {\it vector-valued pluriharmonic
mappings} from $\mathbb{B}^{n}$ into $\mathbb{C}^n$. Let
$f=h+\overline{g}\in\mathcal{PH}(\mathbb{B}^{n}, \mathbb{C}^n)$,
where $h$ and $g$ are holomorphic in $\mathbb{B}^{n}$. Then
$$\det J_{f}=\det\left(\begin{array}{cccc}
\ds Dh & \overline{Dg}  \\
\ds Dg & \overline{Dh}
\end{array}\right)
$$
and if, in addition, $h$ is locally biholomorphic, then one can easily get the formula
$$\det J_{f}=|\det Dh|^{2}\det\left(I_{n}-Dg[Dh]^{-1}\overline{Dg[Dh]^{-1}}\right).
$$
In the case of a {\it planar harmonic mapping} $f=h+\overline{g}$, we find that
$$\det J_f=|f_{z}|^2-|f_{\overline{z}}|^2 
$$
and so, $f$ is locally univalent and sense-preserving in $\ID$ if
and only if $|f_{\overline{z}}(z)|<|f_z(z)|$ in $\ID$; or equivalently if
$f_z(z) \neq0$ and the dilatation $\omega (z)=\overline{f_{\overline{z}}(z)}/f_z(z)$ is analytic
in $\ID$ and has the property that $|\omega (z)|<1$ in $\ID$ (see \cite{Du,Lewy}). For
$f=h+\overline{g}\in\mathcal{PH}(\mathbb{B}^{n}, \mathbb{C}^n)$, the
condition $\|Dg[Dh]^{-1}\|<1$  is sufficient for $\det J_f$ to be
positive and hence for $f$ to be sense-preserving. This is indeed a
natural generalization of one-variable condition.
(cf. \cite{DHK}).

Throughout the discussion a diagonal matrix $A\in
L(\mathbb{C}^{n},\mathbb{C}^{n})$ will be denoted for convenience by
$A=A(\lambda)$ with an understanding that the diagonal entries are
$\lambda _j$, $j=1,2, \ldots, n$, i.e. $(j,j)$-th entry of the
$n\times n $ matrix $A$ is $\lambda _j$. 

\begin{defn}
Let $f=h+\overline{g}$ be a pluriharmonic mapping  from
$\mathbb{B}^{n}$ into $\mathbb{C}^{n}$, where  $h$  and
$g$ are holomorphic in $\mathbb{B}^{n}$. Then
\begin{enumerate}
\item[(1)] $f$ is called  {\it stable pluriharmonic
univalent }  in $\mathbb{B}^{n}$ if
for every   $A\in L(\mathbb{C}^{n},\mathbb{C}^{n})$ with $\|A\|=1$,
the mappings $f_{A}=h+\overline{g}A$ are univalent in $\mathbb{B}^{n}$.
\item[(2)] $f$ is called  {\it stable diagonal
pluriharmonic univalent} in
$\mathbb{B}^{n}$ if for every diagonal matrix $A\in
L(\mathbb{C}^{n},\mathbb{C}^{n})$, the mappings
$f_{A}=h+\overline{g}A$ are univalent in $\mathbb{B}^{n}$ and
$A=A(\lambda)$  with $|\lambda_{j}|=1$ for each $j\in\{1,2,\ldots,n\}$.
\end{enumerate}
The class of all stable pluriharmonic univalent mappings (resp.
stable diagonal pluriharmonic univalent mappings) in
$\mathbb{B}^{n}$ is denoted by {\rm SPU} (resp. {\rm SDPU}) (cf.
\cite{HM}).
\end{defn}

Similarly, we have

\begin{defn}
Let  $F=h+g$ be a holomorphic mapping from $\mathbb{B}^{n}$ into
$\mathbb{C}^{n}$, where $h$ and $g$ are holomorphic in
$\mathbb{B}^{n}$. Then we say that
\begin{enumerate}
\item[(1)] $F$ is {\it stable biholomorphic} in $\mathbb{B}^{n}$ if for every $A\in
L(\mathbb{C}^{n},\mathbb{C}^{n})$ with $\|A\|=1$, the mappings
$h+gA$ are biholomorphic in $\mathbb{B}^{n}$.
\item[(2)] $F$ is {\it stable diagonal biholomorphic}   in $\mathbb{B}^{n}$ if for every  diagonal
matrix $A\in L(\mathbb{C}^{n},\mathbb{C}^{n})$, the mappings $h+gA$
are biholomorphic, where $A=A(\lambda)$ with $|\lambda_{j}|=1$ for each $j\in\{1,2,\ldots,n\}$.
\end{enumerate}
The class of all stable biholomorphic mappings (resp. stable diagonal biholomorphic mappings)
in $\mathbb{B}^{n}$ is denoted by  {\rm SBH} (resp. {\rm SDBH}).
\end{defn}
Now, we state our first result.

\begin{thm}\label{lem-c2}
The pluriharmonic mapping $f=h+\overline{g}$ is in {\rm SPU} {\rm
(}resp. {\rm SDPU)} if and only if the holomorphic mapping $F=h+g$
is in {\rm SBH} {\rm(}resp. {\rm SDBH)}. Moreover, if $f=h+\overline{g}$ is  in {\rm SPU} {\rm (}resp. {\rm
SDPU)} with $\|Dg[Dh]^{-1}\|<1$, where $h$ is locally biholomorphic
and $g$ is holomorphic in $\mathbb{B}^{n}$. Then $h$ is
biholomorphic in $\mathbb{B}^{n}$.
\end{thm}

If we give a  stronger assumption, then we can prove a more general
class of holomorphic mappings $F_{A}=h+gA$ are biholomorphic in
$\mathbb{B}^{n}$, where $A$ is a diagonal matrix with $\|A\|\leq1.$
The result is as follows.

\begin{prop}\label{thm-5}
Let $f=h+\overline{g}$ belong to {\rm SDPU} with
$\|Dg[Dh]^{-1}\|<1$, where $h$ is locally biholomorphic and $g$ is
holomorphic in $\mathbb{B}^{n}$. Then for every $A=A(\lambda)$ with
$|\lambda_{j}|\leq 1$ for each $j\in\{1,2,\ldots,n\}$, the mappings
$F_{A}=h+gA$ are biholomorphic in $\mathbb{B}^{n}$.
\end{prop}

We remark that Proposition \ref{thm-5} is a generalization of
\cite[Theorem 7.1]{HM}.

A domain $D\subset\mathbb{C}^{n}$ is said to be {\it $M$-linearly
connected} if there exists a positive constant $M<\infty$ such that
any two points $w_{1}, w_{2}\in D$ are joined by a path
$\gamma\subset D$ with
$$\ell(\gamma)\leq M\|w_{1}-w_{2}\|,
$$
where $\ell(\gamma)=\inf\left\{\int_{\gamma}\|dz\|:\, \gamma\subset D\right\}.$
It is not difficult to verify that a $1$-linearly connected domain
is convex. For extensive discussions on linearly connected domains,
see \cite{A,CPW,CPW-3,CH,P}.
%

In \cite{CH}, the authors  discussed the relationship between linear
connectivity of the images of $\mathbb{D}$ under the planar harmonic
mappings $f=h+\overline{g}$ and under their corresponding
holomorphic counterparts $h$, where $h$ and $g$ are holomorphic in
$\mathbb{D}$. In \cite[Theorem 5.3]{Clunie-Small-84},  Clunie and
Sheil-Small established an effective and beautiful method of
constructing sense-preserving univalent harmonic mappings defined on
the unit disk, which is popularly called   shear construction. The
following  result is a generalizations of the shearing theorem of
Clunie and Sheil-Small.


\begin{thm}\label{thm-1}
Let $f=h+\overline{g}$ be a  pluriharmonic mapping from
$\mathbb{B}^{n}$ into $\mathbb{C}^{n}$ and $F=h-g$, where  $h$ is
locally biholomorphic  and $g$ is holomorphic in $\mathbb{B}^{n}$.\\

 \item{{\rm (I)}}~ If $F$ is biholomorphic and
$\Omega=F(\mathbb{B}^{n})$ is $M$-linearly connected with
$$\|Dg(z)[Dh(z)]^{-1}\|\leq C<\frac{1}{2M+1} ~\mbox{ for $z\in\mathbb{B}^{n}$},
$$
then for every $A\in L(\mathbb{C}^{n},\mathbb{C}^{n})$ with
$\|A\|\leq1$, the mapping $f_{A}=h+\overline{g}A$ is univalent and
$f_{A}(\mathbb{B}^{n})$ is $M'$-linearly connected, where
$M'=\frac{M(1+C)}{1-(1+2M)C}.$ In particular, $h+g$ is in {\rm
SBH}.\\
\item{{\rm (II)}}~If $f$ is univalent and
$\Omega=f(\mathbb{B}^{n})$ is $M$-linearly connected with
$$\|Dg(z)[Dh(z)]^{-1}\|\leq C<\frac{1}{2M+1} ~\mbox{ for $z\in\mathbb{B}^{n}$},
$$
then for every $A\in L(\mathbb{C}^{n},\mathbb{C}^{n})$ with $\|A\|\leq1$, the mapping
$F_{A}=h-gA$ is univalent and $F_{A}(\mathbb{B}^{n})$ is $M'$-linearly connected,
where $M'=\frac{M(1+C)}{1-(1+2M)C}.$ In particular, $h+\overline{g}$
belongs to {\rm SPU}.
\end{thm}


For $r\in(0,1)$ and a mapping
$f=h+\overline{g}\in\mathcal{PH}(\mathbb{B}^{n}, \mathbb{C}^n)$ with
$\|Dg[Dh]^{-1}\|<1$, the  {\it generalized volume function}
$V_{f}(r)$ of $f$ is defined by
$$V_{f}(r)=\int_{\mathbb{B}^{n}(r)}\|Dh(z)\|^{2n}\left(1-\|Dg(z)[Dh(z)]^{-1}\|^{2}\right)^{n}dV(z),
$$
where $h$ is locally biholomorphic and $g$ is holomorphic in
$\mathbb{B}^{n}$, and $dV$ denotes the normalized Lebesgue volume
measure on $\mathbb{B}^{n}$. We denote by $\mathcal{PH}(V)$, the
class of all mappings
$f=h+\overline{g}\in\mathcal{PH}(\mathbb{B}^{n}, \mathbb{C}^n)$ with
the finiteness condition
$$V:=\sup_{0<r<1}V_{f}(r)<\infty.
$$
We remark that if $n=1$, then
$$V_{f}(r)=\int_{\mathbb{D}_{r}}\det J_{f}(z)\,dA(z)
$$
is the {\it area  function} of the planar harmonic mapping $f$
defined in $\mathbb{D}$, where $dA$ denotes the normalized Lebesgue
area measure on $\mathbb{D}$.

For a holomorphic mapping $f:\,\mathbb{B}^{n}\to \mathbb{C}^{n}$,
$\mathbb{B}^{n}(a,r)$ is called a {\it schlicht ball} of $f$ if
there is a subregion $\Omega\subset\mathbb{B}^{n}$ such that $f$
maps $\Omega$ biholomorphically onto $\mathbb{B}^{n}(a,r)$. We
denote by $B_{f}$ the least upper bound of radii of all schlicht
balls contained in $f(\mathbb{B}^{n})$ and call this the {\it
Landau-Bloch radius} of $f$. The classical theorem of Landau-Bloch
for holomorphic functions in the unit disk fails to extend to
general holomorphic mappings in the ball of $\mathbb{C}^{n}$ (see
\cite{T,W}). However, in 1946, Bochner \cite{B} proved that
Landau-Bloch's theorem does hold for a class of real harmonic
quasiregular mappings. For the extensive discussions on this topic,
see \cite{HG,HG1,CPW-1,CPW6,CPW-2,LX,M3}. Our next aim is to
establish a Landau-Bloch Theorem for $f\in\mathcal{PH}(V)$.

\begin{thm}\label{thm3}
Let $f=h+\overline{g}\in\mathcal{PH}(V)$,
where $V$ is a positive constant, $f(0)=0,$
$Dg(0)=0$ and $|\det J_{f}(0)|=\alpha$ for some  positive constant $\alpha $,
$0<\alpha\leq\frac{8V\psi _{0}^{n}}{\pi} $ with $\ds \psi_{0}=(11+5\sqrt{5})/{2} \approx 11.09017$,
$h$ is locally biholomorphic and $g$ is holomorphic. Then $f$ is
univalent in $\mathbb{B}^{n}(R_u) $, and
$f(\mathbb{B}^{n}(R_u))$ covers the ball
$\mathbb{B}^{n}(R_{c})$, where
$$R_{u} \geq \frac{\alpha\pi (\sqrt{5}-1)(3-2\sqrt{2})}{8V\psi _{0}^{n}}
~\mbox{ and }~ R_{c}\geq\frac{\alpha^{2}\pi(\sqrt{5}-1) (3-2\sqrt{2})}
{16V^{\frac{4n-1}{2n}}\psi_0^{\frac{4n-1}{2}}} 
.
$$
\end{thm}


In fact, the following example will show that there is no
Landau-Bloch theorem for a class of pluriharmonic mappings with
finite volume.

\begin{exam}
For $k\in\{1,2,\ldots\}$ and $z=(z_{1}, \ldots,
z_{n})\in\mathbb{B}^{n}$, let
$$f_{k}(z)=(kz_{1}, z_{2}/k, z_{3},\ldots,z_{n}).
$$
It is not difficult to see that $\det J_{f_{k}}(0)-1=|f_{k}(0)|=0$
and the volume of $f_{k}(\mathbb{B}^{n})$ is
$$\int_{\mathbb{B}^{n}}|\det J_{f_{k}}(z)|\,dV(z)<\infty.
$$
But  there is no absolute constant $s>0$ such that $\mathbb{B}^{n}(s)$ belongs
to $f_{k}(\mathbb{B}^{n})$ for all $k\in\{1,2,\ldots\}$.
\end{exam}

The following result provides a relationship between the real volume
and the generalized volume on a pluriharmonic mapping $f$ defined in $\mathbb{B}^{n}$.

\begin{thm}\label{thm-4}
Let  $f=h+\overline{g}\in\mathcal{PH}(\mathbb{B}^{n}, \mathbb{C}^n)$
with $\|Dg[Dh]^{-1}\|<1$, where $h$ is locally biholomorphic and $g$
is holomorphic in $\mathbb{B}^{n}$. Then for $r\in(0,1)$,
$$\int_{\mathbb{B}^{n}(r)}|\det J_{f}(z)|\,dV(z)\leq K_{r}^{n}V_{f}(r),
$$
where $K_{r}=\frac{1+r}{1-r}$.

\end{thm}

We remark that  Chen and Gauthier proved   Landau-Bloch type Theorems for bounded planar harmonic
mappings (resp. pluriharmonic mappings) (see
\cite[Theorem 3]{HG1} (resp. \cite[Theorem 5]{HG1})). In fact, the volume of all
bounded harmonic mappings (resp. pluriharmonic mappings) is
finite. By Theorem \ref{thm-4}, we see that the condition for functions with finitely generalized  volume
in Theorem \ref{thm3} is weaker than the bounded functions  condition in \cite[Theorems 3 and 5]{HG1}.
In this sense, Theorem \ref{thm3}  generalizes \cite[Theorems 3 and 5]{HG1}.

The proofs of Theorems \ref{thm-5} and \ref{thm-1}  will be
presented in Section \ref{csw-sec2}, and the proof of Theorems
\ref{thm3} and \ref{thm-4}  will be given in Section \ref{csw-sec3}.

\section{univalence criteria on pluriharmonic mappings}\label{csw-sec2}


\subsection*{Proof of Theorem \ref{lem-c2}}
 We begin to prove the necessity of the first part. Let $f=h+\overline{g}$ belong to ${\rm SPU}$. Then, for each
$A\in L(\mathbb{C}^{n},\mathbb{C}^{n})$ with $\|A\|=1$, the mappings
$f_{A}=h+\overline{g}A$ are univalent in $\mathbb{B}^{n}$, where
$h=(h_{1},\ldots,h_{n})$ and $g=(g_{1},\ldots,g_{n})$ are
holomorphic in $\mathbb{B}^{n}$. Suppose on the contrary that
$F=h+g\not\in {\rm SBH}$. Then there exists an $A_{0}\in
L(\mathbb{C}^{n},\mathbb{C}^{n})$ with $\|A_{0}\|=1$ such that
$F_{A_{0}}=h+gA_{0}$ is not biholomorphic in $\mathbb{B}^{n}$. This
means that there are two distinct points
$z_{1},z_{2}\in\mathbb{B}^{n}$ such that \be\label{eq-extra1}
F_{A_{0}}(z_{1})=F_{A_{0}}(z_{2}), ~\mbox{ i.e.
}~h(z_{1})-h(z_{2})=\big(g(z_{2})-g(z_{1})\big)A_{0}. \ee

We divide the rest of the arguments into two cases.

\noindent $\mathbf{ Case ~1.}$ If $h(z_{1})=h(z_{2})$, then
$g(z_{1})=g(z_{2})$. This is a contradiction with the assumption.

\noindent $\mathbf{ Case ~2.}$  If $h(z_{1})\neq h(z_{2})$, then for
$j\in\{1,2,\ldots,n\}$, we may let
$\theta_{j}=\arg\big(h_{j}(z_{1})-h_{j}(z_{2})\big)$ so that
\eqref{eq-extra1} reduces to
$$\big(h(z_{1})-h(z_{2})\big)B=\big(g(z_{2})-g(z_{1})\big)A_{0}B,
$$
which gives
$$\big(h(z_{1})-h(z_{2})\big)B= \overline{\big(h(z_{1})-h(z_{2})\big)}\overline{B}
=\overline{\big(g(z_{2})-g(z_{1})\big)}\overline{A_{0}}\,\overline{B},
$$
or equivalently,
$$ h(z_{1})-h(z_{2}  =\overline{\big(g(z_{2})-g(z_{1})\big)}\overline{A_{0}}\,\overline{B}^2,
$$
where $B=B(\lambda)$ is the diagonal matrix with $\lambda
_j=e^{-i\theta_{j}}$ for $j=1,2,\ldots,n$. This contradicts the
univalency of $f_A=h+\overline{g}A$ for
$A=\overline{A_{0}}\,\overline{B}^{2}$, which proves the necessity
for $f\in{\rm SPU}$. The necessity for $f\in{\rm SDPU}$ is the same,
with $A$ diagonal.

By using similar reasoning as in the proof of the necessity, we can
get the proof of the sufficiency of the first part.

Now we begin to prove the second part of Theorem \ref{lem-c2}. We
prove the second part with a method of contradiction. We suppose
that $h$ is not univalent in $\mathbb{B}^{n}$. Then there are two
distinct points $z_{1},z_{2}\in\mathbb{B}^{n}$ such that
$h(z_{1})=h(z_{2})$. Without loss of generality, we assume that
$z_{1}=h(z_{1})=0.$ In fact, we just take
$$F(z)=\big(h(\phi(z))-h(z_{1})\big)+\overline{\big(g(\phi(z))-g(z_{1})\big)},
$$
where $\phi$ is the automorphism of $\mathbb{B}^{n}$ such that
$\phi(0)=z_{1}$.

Therefore, we can assume that $f=h+\overline{g}$ is in {\rm SPU}
{\rm (}resp. {\rm SDPU)} with the normalized conditions
$h(0)=g(0)=0.$ By the conclusion of the first part of Theorem
\ref{lem-c2}, we obtain that $h+g$ belongs to {\rm SBH} and
therefore, it follows that for each
$z\in\mathbb{B}^{n}\backslash\{0\}$ either $\|h(z)\|<\|g(z)\|$ or
$\|g(z)\|<\|h(z)\|$.
Now, $h(z_2)=0$ so  $\|h(z_2)\|<\|g(z_2)\|$. Then by the continuity of the function
$\|h(z)\|-\|g(z)\|$, we conclude that $\|h(z)\|<\|g(z)\|$ in $\mathbb{B}^{n}\backslash\{0\}$.

Since $h$ is locally biholomorphic, we see that there is a sequence of points in $\{Z_{n}\}_{n\geq 1}$ in
$\mathbb{B}^{n}\backslash\{0\}$ such that $Z_{n}\rightarrow 0$  as $n\rightarrow\infty.$
Thus,
$$\frac{\|h(Z_{n})\|}{\|Z_{n}\|} = \left\| \frac{h(Z_{n})-h(0)}{\|Z_{n}-0\|}\right\| < \left\| \frac{g(Z_{n})-g(0)}{\|Z_{n}-0\|}\right\|
=\frac{\|g(Z_{n})\|}{\|Z_{n}\|} ~\mbox{ for $n\geq 1$,}
$$
and therefore, we obtain that
\be\label{eq-t}
\|Dh(0)\|\leq\|Dg(0)\|.
\ee
On the other hand, by the assumption, we have
$$\|Dg(0)\|=\|Dg(0)[Dh(0)]^{-1}Dh(0)\|\leq\|Dg(0)[Dh(0)]^{-1}\|\|Dh(0)\|<\|Dh(0)\|
$$
which contradicts the inequality (\ref{eq-t}). The proof of this
theorem is complete. \qed

\subsection*{Proof of Proposition \ref{thm-5}}
By the assumption of Proposition \ref{thm-5} and Theorem
\ref{lem-c2}, we obtain that the mappings $F_{A}=h+gA$ are
biholomorphic for every diagonal matrix $A\in
L(\mathbb{C}^{n},\mathbb{C}^{n})$, where $A=A(\lambda)$ with
$|\lambda_{j}|=1$ for each $j=1,2,\ldots,n$.

Next we prove that  for every diagonal matrix $A\in
L(\mathbb{C}^{n},\mathbb{C}^{n})$, where $A=A(\lambda)$ with
$|\lambda_{j}|<1$ for some $j\in\{1,2,\ldots,n\}$, the mappings
$F_{A}=h+gA$ are also biholomorphic.  Without loss of generality, we
assume that diagonal matrix $A(\lambda)\in
L(\mathbb{C}^{n},\mathbb{C}^{n})$ with $|\lambda_{j}|<1$ for each
$j=1,2,\ldots,n$. In other words, we will prove that if
$f=h+\overline{g}$ belongs to {\rm SDPU} and the diagonal matrix
$A=A(\lambda)$ is with $|\lambda_{j}|<1$ for each $j=1,2,\ldots,n$,
then we have $T_{A}(f)=f+\overline{f}A \in {\rm SDPU}$.

By simple calculations, we have
$$T_{A}(f)=h+gA+\overline{g+h\overline{A}}.
$$
By Theorem \ref{lem-c2}, in order to prove
$T_{A}(f)=f+\overline{f}A$ is in {\rm SDPU}, we only need to prove
that for every diagonal matrix of the form $D=D(\lambda )$ with
diagonal entries $\lambda _j=e^{i\theta _j}$ ($j=1,2,\ldots,n$),
the mappings $h+gA+(g+h\overline{A})D$ are biholomorphic in $\mathbb{B}^{n}$,
where  $\theta_{j}\in[0,2\pi)$ for each $j=1,2,\ldots,n$. We may rewrite
\beqq
h+gA+(g+h\overline{A})D &= &g(A+D)+h(I_{n}+\overline{A}D)\\
&=&\big[h+g(A+D)(I_{n}+\overline{A}D)^{-1}\big](I_{n}+\overline{A}D).
\eeqq
Then it is a simple exercise to see that the matrix $C=(A+D)(I_{n}+\overline{A}D)^{-1}$ is a diagonal matric with diagonal entries
$$ \varphi_j(\lambda_j)= \frac{\lambda_{j}+e^{i\theta_{j}}}{1+\overline{\lambda}_{j}e^{i\theta_{j}}} ~\mbox{ for  $j=1,2,\ldots,n.$
}
$$
Since $|\varphi_j(\lambda_j)|=1$ for each $j=1,2,\ldots,n$, we
conclude that
$$\|C\|=\|(A+D)(I_{n}+\overline{A}D)^{-1}\|=1
$$
and that the mapping $h+gA+(g+h\overline{A})D$  belongs to {\rm
SDBH} so that $T_{A}(f)=f+\overline{f}A$ is in {\rm SDPU}. It
follows from Theorem \ref{lem-c2} that the holomorphic part $h+gA$
of $T_{A}(f)$ is biholomorphic. The proof of this proposition is
complete. \qed

\vspace{8pt}


The following lemma plays a key role in the proofs of Theorem
\ref{thm-1}. 

\begin{Lem}\label{lem-1}
Let $A$ be an $n\times n$ complex matrix with $\|A\|<1$. Then
$I_{n}\pm A$ are nonsingular matrixes and $\|(I_{n}\pm A)^{-1}\|\leq1/{(1-\|A\|)}.$
\end{Lem}

Lemma \Ref{lem-1} may be referred to as the Neumann expansion
theorem,  and so the proof is omitted here.

\subsection*{Proof of Theorem \ref{thm-1}}

We first prove part {\rm (I)}.
 Let $A\in L(\mathbb{C}^{n},\mathbb{C}^{n})$ with $\|A\|\leq1$ and
consider the mapping $f_{A}=h+\overline{g}A$. Using the hypotheses,
we first show that $f_{A}(\mathbb{B}^{n})$ is linearly connected.
Define $\Omega=F(\mathbb{B}^{n})$, where $F=h-g$,  and
$$H(w)=f_{A}(F^{-1}(w))=w+g(F^{-1}(w))+\overline{g(F^{-1}(w))}A ~\mbox{ for $w\in\Omega.$}
$$
Clearly
$$DF=Dh-Dg=(I_{n}-Dg[Dh]^{-1})Dh,
$$
and therefore, we see that
$$[DF]^{-1}=[Dh]^{-1}\left(I_{n}-Dg[Dh]^{-1}\right)^{-1}
$$
which gives \be\label{eqh1}
DH=I_{n}+Dg[DF]^{-1}=I_{n}+Dg[Dh]^{-1}\left(I_{n}-Dg[Dh]^{-1}\right)^{-1}
\ee and \be\label{eqh2}
\overline{D}H=\overline{Dg[DF]^{-1}}A=\overline{Dg[Dh]^{-1}}\overline{\left(I_{n}-Dg[Dh]^{-1}\right)^{-1}}A.
\ee

 For any two distinct points
$w_{1},\, w_{2}\in\Omega$, by hypothesis, there is a curve
$\gamma\subset\Omega$ joining $w_{1}$ and $w_{2}$ such that
$l(\gamma)\leq M\|w_{1}-w_{2}\|$. Also, we let $\Gamma=H(\gamma)$.
On one hand, by (\ref{eqh1}), (\ref{eqh2}) and Lemma \Ref{lem-1}, we find
that
\beq\label{eq2-extra}
\nonumber
l(\Gamma)&=&\int_{\Gamma}\|dH(w)\|\leq\int_{\gamma}\big(\|DH(w)\|+\|\overline{D}H(w)\|\big)\,\|dw\|\\
\nonumber &\leq&\int_{\gamma}\big(\|I_{n}\|+2\big\|Dg[Dh]^{-1}\left(I_{n}-Dg[Dh]^{-1}\right)^{-1}\big\|\big)\,\|dw\|\\
\nonumber
&\leq&\int_{\gamma}\left(1+\frac{2\|Dg[Dh]^{-1}\|}{1-\|Dg[Dh]^{-1}\|}\right)\,\|dw\|\\
&\leq&\frac{1+C}{1-C}M\|w_{2}-w_{1}\|.
\eeq
On the other hand, the definition of $H$ gives
\beq\label{eq-1-th3}
\|H(w_{2})-H(w_{1})\|&\geq&\|w_{2}-w_{1}\|-2\|g(F^{-1}(w_{2}))-g(F^{-1}(w_{1}))\|\\
\nonumber &\geq&\|w_{2}-w_{1}\|-2\int_{\gamma}\|Dg[DF]^{-1}\|\,\|dw\|\\
\nonumber &\geq&\|w_{2}-w_{1}\|-2\int_{\gamma} \frac{\|Dg[Dh]^{-1}\|}{1-\|Dg[Dh]^{-1}\|} \, \|dw\|\\
\nonumber &\geq&\frac{1-C(1+2M)}{1-C}\|w_{2}-w_{1}\|,
\eeq
and therefore, \eqref{eq2-extra} gives
$$l(\Gamma)\leq M'\|H(w_{2})-H(w_{1})\|,
$$
where $M'=\frac{(1+C)M}{1-(1+2M)C}.$

Finally,  we show that $f_{A}=h+\overline{g}A$ is univalent for
every $A\in L(\mathbb{C}^{n},\mathbb{C}^{n})$ with $\|A\|\leq1$.
Suppose on the contrary that $f_{A}$ fails to be univalent. Then there are two
distinct points $w_{1},w_{2}$ such that $H(w_{1})=H(w_{2})$ which
is impossible, by (\ref{eq-1-th3}).

At last, since $f_{A}=h+\overline{g}A$ is univalent for every $A\in
L(\mathbb{C}^{n},\mathbb{C}^{n})$ with $\|A\|\leq1$, it follows from
Lemma \ref{lem-c2} that $h+g$ belongs to {\rm SBH} and the proof of
 part {\rm (I)} is finished.

Now we begin to prove  part {\rm (II)}. Assume the hypotheses and
Define
$$H(w)=F_{A}(f^{-1}(w))=w-\overline{G(w)}-G(w)A,
$$
where $G=g\circ f^{-1}$ and $w=f(z)$. By the chain rule, we have
$$DG=Dg Df^{-1}~\mbox{ and }~\overline{D}G=Dg \overline{D}f^{-1},
$$
which implies
$$DH=I_{n}-\overline{Dg}~ \overline{\overline{D}f^{-1}}-Dg Df^{-1}A~\mbox{ and }~
\overline{D}H=-\overline{Dg}~ \overline{Df^{-1}}-Dg
\overline{D}f^{-1}A.
$$
It follows from the inverse mapping theorem and Lemma \Ref{lem-1}
that $f^{-1}$ is differentiable. Differentiation of the equation
$f^{-1}(f(z))=z
$
yields the following relations
$$\begin{cases}
\displaystyle Df^{-1} Dh+\overline{D}f^{-1} Dg=I_{n},\\
\displaystyle Df^{-1} \overline{Dg}+\overline{D}f^{-1}
\overline{Dh}=0,
\end{cases}
$$
which give
\be\label{eq-2}
\begin{cases}
\displaystyle Df^{-1}=[Dh]^{-1}\left(I_{n}-\overline{Dg}[\overline{Dh}]^{-1}Dg[Dh]^{-1}\right)^{-1},\\
\displaystyle
\overline{D}f^{-1}=-[Dh]^{-1}\left(I_{n}-\overline{Dg}[\overline{Dh}]^{-1}Dg[Dh]^{-1}\right)^{-1}\overline{Dg}[\overline{Dh}]^{-1}.
\end{cases}
\ee
By (\ref{eq-2}) and Lemma \Ref{lem-1}, we get
\begin{eqnarray*}
\|DG\|+\|\overline{D}G\|&=& \|DgDf^{-1}\|+\|Dg\overline{D}f^{-1}\|\\
&=&\big\|Dg[Dh]^{-1}\left(I_{n}-\overline{Dg}[\overline{Dh}]^{-1}Dg[Dh]^{-1}\right)^{-1}\big\|\\
&&+\big\|Dg[Dh]^{-1}\left(I_{n}-\overline{Dg}[\overline{Dh}]^{-1}Dg[Dh]^{-1}\right)^{-1}\overline{Dg}[\overline{Dh}]^{-1}\big\|\\
&\leq&\big\|\left(I_{n}-\overline{Dg}[\overline{Dh}]^{-1}Dg[Dh]^{-1}\right)^{-1}\big\|\|Dg[Dh]^{-1}\|\\
&&\times \left(1+\|Dg[Dh]^{-1}\|\right)\\
&\leq&\frac{\|Dg[Dh]^{-1}\|\left(1+\|Dg[Dh]^{-1}\|\right)}
{1-\big\|\overline{Dg}[\overline{Dh}]^{-1}Dg[Dh]^{-1}\big\|}\\
&\leq&\frac{\|Dg[Dh]^{-1}\|\left(1+\|Dg[Dh]^{-1}\|\right)}{1-\|Dg[Dh]^{-1}\|^{2}}\\
&\leq&\frac{\|Dg[Dh]^{-1}\|}{1-\|Dg[Dh]^{-1}\|}\\
&<&\frac{C}{1-C}.
\end{eqnarray*}
Let $\gamma\subset\Omega$ be a curve joining $w_{1}, w_{2}$ with $l(\gamma)\leq M\|w_{1}-w_{2}\|$.

Then
\begin{eqnarray*}
\|DH\|+\|\overline{D}H\|&=&\big\|I_{n}-\overline{Dg}~
\overline{\overline{D}f^{-1}}-Dg Df^{-1}A\big\|+\big\|\overline{Dg}~
\overline{Df^{-1}}+Dg
\overline{D}f^{-1}A\big\|\\
&\leq&1+2\big\|\overline{Dg}~
\overline{\overline{D}f^{-1}}\big\|+2\big\|\overline{Dg}~
\overline{Df^{-1}}\big\|\\
&\leq&1+\frac{2\|Dg[Dh]^{-1}\|}{1-\|Dg[Dh]^{-1}\|}\\
&=&\frac{1+\|Dg[Dh]^{-1}\|}{1-\|Dg[Dh]^{-1}\|}\\
&\leq&\frac{1+C}{1-C},
\end{eqnarray*}
which gives
\beq\label{eq1-extra}
\nonumber l\big(H(\gamma)\big)&\leq&\int_{\gamma}\big(\|DH(w)\|+\|\overline{D}H(w)\|\big)\,\|dw\|\\
\nonumber &\leq& \frac{1+C}{1-C}l(\gamma)\\
&\leq&\frac{M(1+C)}{1-C}\|w_{2}-w_{1}\|.
\eeq

On the other hand, by \eqref{eq1-extra},
\beq\label{eq7h} \nonumber
\|H(w_{2})-H(w_{1})\|& =&\|w_{2}-\overline{G(w_{2})}-G(w_{2})A -\big(w_{1}-\overline{G(w_{1})}-G(w_{1})A\big) \|\\
\nonumber&\geq&\|w_{2}-w_{1}\|-2\|G(w_{2})-G(w_{1})\|\\ \nonumber
&\geq&\|w_{2}-w_{1}\|-2\int_{\gamma}\big(\|DG(w)\|+\|\overline{D}G(w)\|\big)\,\|dw\|\\\nonumber
&\geq&\|w_{2}-w_{1}\|-\frac{2C}{1-C}l(\gamma)\\
&\geq&\frac{1-C(1+2M)}{1-C}\|w_{2}-w_{1}\|\\ \nonumber
&\geq&[1-C(1+2M)]\frac{l\big(H(\gamma)\big)}{M(1+C)}.
\eeq
Hence
$$l\big(H(\gamma)\big)\leq\frac{M(1+C)}{1-C(1+2M)}\|H(w_{2})-H(w_{1})\|.
$$
Now we prove the univalency of $F_{A}$ for every $A\in
L(\mathbb{C}^{n},\mathbb{C}^{n})$ with $\|A\|\leq1$.  Suppose that
$H(w_{1})=H(w_{2})$ for distinct points $w_{1},w_{2}\in
f(\mathbb{B}^{n})$. Then there exist two distinct points $w_{1}$,
$w_{2}$ such that $H(w_{1})=H(w_{2})$ which is a contradiction to
(\ref{eq7h}).

At last, since $F_{A}$ is biholomorphic for every $A\in
L(\mathbb{C}^{n},\mathbb{C}^{n})$ with $\|A\|\leq1$, it follows
Lemma \ref{lem-c2} that $h+\overline{g}$ is in {\rm SPU}. The proof
of the theorem is complete. \qed


\section{The Landau-Bloch theorem on pluriharmonic mappings}\label{csw-sec3}
For $f \in \mathcal{H}(\IB^n, \mathbb{C}^n)$, we use the
standard notations:
\be\label{eqbe1}
\Lambda_{f}(z)=\max_{ \theta\in\partial\mathbb{B}^{n}}\|Df(z)\theta+\overline{D}f(z)\overline{\theta}\|\;\;
\mbox{ and }\;\; \lambda_{f}(z)=\min_{
\theta\in\partial\mathbb{B}^{n}}\|Df(z)\theta+\overline{D}f(z)\overline{\theta}\|.
\ee
We see that (see for instance \cite{HG1,CPW-2})
\be\label{eqbe2}
\Lambda_{f}=\max_{\theta\in\partial\mathbb{B}^{2n}_{\mathbb{R}}}\|J_{f}\theta\|\;\;
\mbox{ and }\;\;
\lambda_{f}=\min_{\theta\in\partial\mathbb{B}^{2n}_{\mathbb{R}}}\|J_{f}\theta\|.
\ee
Then the following two results are useful for the proof of Theorem \ref{thm3}.

\begin{Thm}{\rm (\cite[Theorem 4]{HG1})}\label{ThmB}
Let $f$ be a pluriharmonic mappings of $\mathbb{B}^{n}$ into
$\mathbb{B}^{m}$. Then
$$\Lambda_{f}(z)\leq\frac{4}{\pi}\frac{1}{1-\|z\|^{2}}~\mbox{ for }~z\in\mathbb{B}^{n}.
$$
If $f(0)=0,$ then
$\|f(z)\|\leq (4/\pi)\arctan\|z\|$  for $z\in\mathbb{B}^{n}.$
\end{Thm}

\begin{Lem} {\rm (\cite[Lemma 1]{HG1} or \cite[Lemma 4]{LX})}\label{LemA}
Let $A$ be an $n \times n$ complex $($real$)$ matrix with
$\|A\|\neq0$. Then for any unit vector $\theta\in\partial
\mathbb{B}^{n}$, the inequality $\|A\theta\|\geq |\det
A|/\|A\|^{n-1} $ holds.
\end{Lem}

\subsection*{Proof of Theorem \ref{thm3}}
For each fixed $\theta\in\partial\mathbb{B}^{n}$, let
$A_{\theta}=Dg[Dh]^{-1}\theta$. Since $\|Dg[Dh]^{-1}\|<1$,   by
Schwarz's lemma, we see that for $z\in\mathbb{B}^{n}(r)$,
$\|A_{\theta}(z)\|<\|z\|$ if $r\in(0,1)$. The arbitrariness of
$\theta\in\partial\mathbb{B}^{n}$ gives
\be\label{eq-3}
\|Dg(z)[Dh(z)]^{-1}\|<r
\ee
for $z\in\mathbb{B}^{n}(r)$.

By (\ref{eq-3}), we get
\be\label{eq-3.1}
\frac{1+\|Dg(z)[Dh(z)]^{-1}\|}{1-\|Dg(z)[Dh(z)]^{-1}\|}\leq\frac{1+r}{1-r}=K_{r},
\quad  z\in\mathbb{B}^{n}(r).
\ee
Applying (\ref{eq-3.1}), we obtain that
\begin{eqnarray*}
V_{f}(r)&=&
\int_{\mathbb{B}^{n}(r)}\|Dh(\zeta)\|^{2n}\big(1-\|Dg(z)[Dh(z)]^{-1}\|^{2}\big)^{n}\,dV(\zeta)\\
&\geq&\frac{1}{K_{r}^{n}}\int_{\mathbb{B}^{n}(r)}\|Dh(\zeta)\|^{2n}\big(1+\|Dg(z)[Dh(z)]^{-1}\|\big)^{2n}\,dV(\zeta).
\end{eqnarray*}
Fix $z\in\mathbb{B}^{n}$ and let
$D_{z}^{r}=\{\zeta\in\mathbb{C}^{n}:\, \|\zeta-z\|<r-\|z\|\}.$  Since
$\|Dh\theta+\overline{Dg}\,\overline{\theta}\|^{2n}$ is subharmonic
in $\mathbb{B}^{n}$, we see that for $\rho\in[0,1-\|z\|)$,
\be\label{eq-4}
\|Dh(z)\theta+\overline{Dg(z)}\,\overline{\theta}\|^{2n}\leq\int_{\partial\mathbb{B}^{n}}
\|Dh(z+\rho\zeta)\theta+\overline{Dg(z+\rho\zeta)}\, \overline{\theta}\|^{2n}\,d\sigma(\zeta),
\ee
where $\theta\in\partial\mathbb{B}^{n}$. Multiplying both sides of (\ref{eq-4}) by $2nr^{2n-1}$ and then
integrating from $0$ to $r-\|z\|$, we obtain

\vspace{7pt}
$\ds (r-\|z\|)^{2n}\|Dh(z)\theta+\overline{Dg(z)}\, \overline{\theta}\|^{2n}$
\begin{eqnarray*}
&\leq&\int_{0}^{r-\|z\|}\left[2n\rho^{2n-1}
\int_{\partial\mathbb{B}^{n}}\|Dh(z+\rho\zeta)\theta+\overline{Dg(z+\rho\zeta)}\,\overline{\theta}\|^{2n}d\sigma(\zeta)\right]d\rho\\
&=&\int_{D_{z}^{r}}\|Dh(\zeta)\theta+\overline{Dg(\zeta)}\, \overline{\theta}\|^{2n}\,dV(\zeta)\\
&\leq&\int_{\mathbb{B}^{n}(r)}\|Dh(\zeta)\theta+\overline{Dg(\zeta)}\, \overline{\theta}\|^{2n}\,dV(\zeta)
\\
&\leq&\int_{\mathbb{B}^{n}(r)}(\|Dh(\zeta)\|+\|\overline{Dg(\zeta)}\|)^{2n}\,dV(\zeta)\\
&\leq&\int_{\mathbb{B}^{n}(r)}(1+\|Dg(\zeta)[Dh(\zeta)]^{-1}\|)^{2n}\|Dh(\zeta)\|^{2n}\,dV(\zeta)\\
&\leq&K_{r}^{n}V_{f}(r)\\
&\leq&K_{r}^{n}V,
\end{eqnarray*}
which implies that
\be\label{eq-3.2}
\Lambda_{f}^{2n}(z)=\max_{\theta\in\partial\mathbb{B}^{n}}\|Dh(z)\theta+\overline{Dg(z)}\, \overline{\theta}\|^{2n}
\leq\frac{K_{r}^{n}V}{(r-|z|)^{2n}}.
\ee
For $\xi\in\mathbb{B}^{n}$, let $F(\xi)=r^{-1}f(r\xi)=H(\xi)+\overline{G(\xi)}$,
where
$$H(\xi)=r^{-1}h(r\xi) ~\mbox{ and }~ G(\xi)=r^{-1}g(r\xi).
$$
Then
\be\label{eq-3.3}
\Lambda_{F}(\xi)=\max_{\theta\in\partial\mathbb{B}^{n}}\|DH(\xi)\theta+\overline{DG(\xi)}\, \overline{\theta}\|
\leq\frac{K_{r}^{\frac{1}{2}}V^{\frac{1}{2n}}}{r(1-|\xi|)},
\ee
where $z=r\xi.$ By (\ref{eq-3.3}), we have
\be\label{eq-3.30}
\Lambda_{F}(\xi)\leq\frac{V^{\frac{1}{2n}}\sqrt{\min_{0<r<1}\psi(r)}}{1-|\xi|}
=\frac{V^{\frac{1}{2n}}\sqrt{\psi(r_{0})}}{1-|\xi|},
\ee
where
$$\psi(r)=\frac{1+r}{r^{2}(1-r)}~\mbox{ and }~ r_{0}=\frac{\sqrt{5}-1}{2}.
$$
Again, for $w\in\mathbb{B}^{n}$ and a fixed $t\in(0,1)$, let
$P(w)=t^{-1}F(tw).$ Then
\be\label{eq-3.5}
\Lambda_{P}(w)=\Lambda_{F}(\xi)\leq
\frac{V^{\frac{1}{2n}}\sqrt{\psi(r_{0})}}{1-|\xi|}=\frac{V^{\frac{1}{2n}}\sqrt{\psi_{0}}}{1-t|w|}
\leq\frac{V^{\frac{1}{2n}}\sqrt{\psi_{0}}}{1-t}=M(t)
\ee
where $\psi_0=\psi(r_{0})$ and by a computation, it follows easily that
$$ \psi_{0}=\frac{11+5\sqrt{5}}{2}   \approx 11.09017.
$$
Let $w_{1},$ $w_{2}$ be two distinct points in $\mathbb{B}^{n}(\rho(t))$
with $w_{2}-w_{1}=\|w_{2}-w_{1}\|\theta$, where
$$\rho(t)=\frac{\alpha\pi}{4(M(t)+M(0))(M(0))^{2n-1}}=\frac{\alpha\pi(1-t)}{4(M(0))^{2n}(2-t)}=\frac{\alpha\pi(1-t)}{4V\psi_0^{2n}(2-t)}.
$$
By the assumption and monotonicity of $\rho(t)$, it can be easily see that for
 $\rho(t)\leq\rho(0)\leq1$ for $t\in[0,1)$. Define the pluriharmonic
mapping
$$\phi_{\theta}(w)=(DP(w)-DP(0))\theta+(\overline{D}P(w)-\overline{D}P(0))\overline{\theta} ~\mbox{ for }~ w\in\mathbb{B}^{n}.
$$
By (\ref{eq-3.5}), we get
$$\|\phi_{\theta}(w)\|\leq\Lambda_{P}(w)+\Lambda_{P}(0)\leq M(t)+M(0)
$$
and therefore, by using Theorem \Ref{ThmB},  we obtain
\be\label{eq3.5}
\|\phi_{\theta}(w)\|\leq\frac{4(M(t)+M(0))}{\pi}\arctan\|w\|\leq\frac{4(M(t)+M(0))}{\pi}\|w\|
\ee
for  $w\in\mathbb{B}^{n}$. By (\ref{eqbe1}), (\ref{eqbe2}), (\ref{eq-3.30}) and Lemma \Ref{LemA}, we have
\be\label{eq-3.6}
\lambda_{P}(0)\geq\frac{|\det
J_{P}(0)|}{\Lambda_{P}^{2n-1}(0)} \geq\frac{\alpha}{(M(0))^{2n-1}}.
\ee

Let $[w_{1},w_{2}]$ denote the segment from $w_{1}$ to $w_{2}$,
$dw=(du_{1},\ldots,du_{n})^{T}$ and
$d\overline{w}=(d\overline{u}_{1},\ldots,d\overline{u}_{n})^{T}$.
Then by (\ref{eq3.5}) and (\ref{eq-3.6}) we have
\begin{eqnarray*}
\|P(w_{1})-P(w_{2})\|&=&\left\|\int_{[w_{1},w_{2}]}DP(w)\,dw+
\overline{D}P(w)\,d\overline{w}\right\|\\
&\geq&\left\|\int_{[w_{1},w_{2}]}DP(0)\,dw+
\overline{D}P(0)\,d\overline{w}\right\|\\
&&-\left\|\int_{[w_{1},w_{2}]}(DP(w)-DP(0))\,dw+
(\overline{D}P(w)-\overline{D}P(0))\,d\overline{w}\right\|\\
&\geq&\|w_{1}-w_{2}\|\lambda_{P}(0)-\int_{[w_{1},w_{2}]}\|\phi_{\theta}(w)\|\,\|dw\|\\
&>&\|w_{1}-w_{2}\|\left[\frac{\alpha}{(M(0))^{2n-1}}-\frac{4(M(t)+M(0))}{\pi}\rho(t)\right]=0.
\end{eqnarray*}

Furthermore, for $w$ with $|w|=\rho(t)$, we have
\begin{eqnarray*}
\|P(w)-P(0)\|&\geq&
\left\|\int_{[0,w]}DP(0)d\varsigma+\overline{D}P(0)\,d\overline{\varsigma}\right\|\\
&&
-\left\|\int_{[0,w]}(DP(\varsigma)-DP(0))\,d\varsigma+(\overline{D}P(\varsigma)-
(\overline{D}P(0))\,d\overline{\varsigma}\right\|\\
&>&\rho(t)\left[\frac{\alpha}{(M(0))^{2n-1}}-\frac{2(M(t)+M(0))}{\pi}\rho(t)\right]\\
&=&\frac{\alpha^{2}\pi}{8(M(t)+M(0))(M(0))^{4n-2}}\\
&=&\frac{\alpha^{2}\pi(1-t)}{8(M(0))^{4n-1}(2-t)}.
\end{eqnarray*}
The last inequality shows that the range
$P(\mathbb{B}^{n}(\rho(t)))$ contains a univalent ball
$\mathbb{B}^{n}(R)$, where
$$R=R(t)\geq\frac{\alpha^{2}\pi(1-t)}{8(M(0))^{4n-1}(2-t)}.
$$
Therefore, $f$ is univalent in $\mathbb{B}^{n}(tr_{0}\rho(t))$, and
$f(\mathbb{B}^{n}(tr_{0}\rho(t)))$  covers the ball
$\mathbb{B}^{n}(R_{c})$ with
\beqq
R_{c}= \max_{0<t<1}(r_{0}tR(t))&\geq&\frac{\alpha^{2}\pi
r_{0}\max_{0<t<1} \nu(t)}{8(M(0))^{4n-1}}\\
&=&\frac{\alpha^{2}\pi r_{0}(3-2\sqrt{2})}{8(M(0))^{4n-1}} \\
&= &\frac{\alpha^{2}\pi (\sqrt{5} -1) (3-2\sqrt{2})}
{16V^{\frac{4n-1}{2n}}\psi_0^{\frac{4n-1}{2}}},
\eeqq
where $\nu(t)=\frac{t(1-t)}{(2-t)}$ and a computation gives that
$$\max_{0<t<1}\nu(t)=\nu(2-\sqrt{2})=3-2\sqrt{2}.
$$
Moreover, using the value of $r_0$ and the function $\rho (t)$, we find that
$$tr_{0}\rho(t) =\frac{\alpha\pi (\sqrt{5}-1)}{8V\psi _{0}^{n}}\nu(t)
\geq \frac{\alpha\pi (\sqrt{5}-1)(3-2\sqrt{2})}{8V\psi _{0}^{n}} =R_u.
$$
The proof of this theorem is complete.\qed

\subsection*{Proof of Theorem \ref{thm-4}}
Applying Lemma \Ref{LemA}, for an $n \times n$ complex matrix A, we
have
\be\label{m-1} |\det A|\leq\|A\|^{n}.
\ee
By (\ref{eq-3.1}) and (\ref{m-1}), for $z\in\mathbb{B}^{n}(r)$, we get
\begin{eqnarray*}
|\det J_{f}(z)|&=&|\det
Dh(z)|^{2}\left|\det\big(I_{n}-Dg(z)[Dh(z)]^{-1}\overline{Dg(z)[Dh(z)]^{-1}}\big)\right|\\
&\leq&\|Dh(z)\|^{2n}\left\|\big(I_{n}-Dg(z)[Dh(z)]^{-1}\overline{Dg(z)[Dh(z)]^{-1}}\big)\right\|^{n}\\
&\leq&\|Dh(z)\|^{2n}\left(1+\|Dg(z)[Dh(z)]^{-1}\|^{2}\right)^{n}\\
&\leq&\|Dh(z)\|^{2n}\left(1+\|Dg(z)[Dh(z)]^{-1}\|\right)^{2n}\\
&\leq&\left(\frac{1+r}{1-r}\right)^{n}\|Dh(z)\|^{2n}\left(1-\|Dg(z)[Dh(z)]^{-1}\|^{2}\right)^{n},
\end{eqnarray*}
which implies that
\begin{eqnarray*}
\int_{\mathbb{B}^{n}(r)}|\det J_{f}(z)|dV(z)&\leq&K_{r}^{n}\int_{\mathbb{B}^{n}(r)}
\|Dh(z)\|^{2n}\left(1-\|Dg(z)[Dh(z)]^{-1}\|^{2}\right)^{n}dV(z)\\
&=&K_{r}^{n}V_{f}(r),
\end{eqnarray*}
where $K_{r}=\frac{1+r}{1-r}$. The proof of this theorem is
complete.\qed



\normalsize


\begin{thebibliography}{99}

\bibitem{A} {\sc J. M. Anderson, J. Becker and  J. Gevirtz,}
First-order univalence criteria interior chord-arc conditions, and quasidisks,
\textit{Michigan  Math. J.}, {\bf 56}(2008), 623--636.

\bibitem{B}  {\sc S.~Bochner,}  Bloch's theorem for real variables,
\textit{Bull. Amer. Math. Soc., } {\bf 52}(1946),  715--719.

\bibitem{HG}
 {\sc H.  Chen and P. M. Gauthier,}
 Bloch constants in several variables,
\textit{Trans. Amer. Math. Soc.,} {\bf 353} (2001), 1371--1386.

\bibitem{HG1} {\sc H. Chen and P. M. Gauthier,}
The Landau theorem and Bloch theorem for planar harmonic and
pluriharmonic mappings,
\textit{Proc. Amer. Math. Soc.,} {\bf 139}(2011), 583--595.

\bibitem{CPW} {\sc Sh.~Chen, S.~Ponnusamy and X.~Wang,}
Properties of some classes of planar harmonic and planar biharmonic mappings,
\textit{Complex Anal. Oper. Theory,} {\bf 5}(2011), 901--916.

\bibitem{CPW-1} {\sc Sh. Chen,  S. Ponnusamy and X. Wang,}
Equivalent moduli of continuity, Bloch's theorem for pluriharmonic mappings in $\mathbb{B}^{n}$,
\textit{Proc. Indian Acad. Sci. {\rm (}Math. Sci.{\rm)}}, {\bf 122}(2012), 583--595.

\bibitem{CPW6} {\sc Sh. Chen,  S. Ponnusamy and X. Wang,}
Weighted Lipschitz continuity, Schwarz-Pick's Lemma and
Landau-Bloch's theorem for hyperbolic harmonic mappings in $\mathbb{C}^{n}$,
\textit{Math. Model. Anal.}, {\bf 18}(2013), 66--79.

\bibitem{CPW-2} {\sc Sh.~ Chen, S.~Ponnusamy and X.~Wang,}
Harmonic mappings in Bergman spaces,
\textit{Monatsh. Math.}, {\bf 170}(2013), 325--342.


\bibitem{CPW-3} {\sc Sh.~ Chen, S.~Ponnusamy and X.~Wang,}
Univalence criteria and Lipschitz-type spaces on pluriharmonic mappings,
\textit{ Math. Scand.}, to appear.



\bibitem{CH} {\sc M. Chuaqui and R. Hern\'andez,}
Univalent  harmonic mappings and linearly connected domains,
\textit{J. Math. Anal. Appl.,} {\bf 332}(2007), 1189--1194.

\bibitem{Clunie-Small-84} {\sc J. G.  Clunie and T. Sheil-Small,}
Harmonic univalent functions,
\textit{Ann. Acad. Sci. Fenn. Ser. A I Math.,} {\bf 9}(1984), 3--25.



\bibitem{Du} {\sc P. Duren,}
{\it Harmonic mappings in the plane,} Cambridge Univ. Press, 2004.

\bibitem{DHK} {\sc P. Duren, H. Hamada and G. Kohr,}
Two-point distortion theorems for harmonic and pluriharmonic mappings,
\textit{Trans. Amer. Math. Soc.}, {\bf 363}(2011), 6197--6218.

\bibitem{HM} {\sc  R. Hern\'andez and M. J. Mart\'in,}
Stable geometric properties of analytic and harmonic functions,
\textit{Proc. Cambridge Phil. Soc.},  {\bf 155}(2013), 343--359.

\bibitem{I} {\sc A. J. Izzo,}
Uniform algebras generated by holomorphic and pluriharmonic functions,
\textit{Trans. Amer. Math. Soc.}, {\bf 339}(1993), 835--847.


\bibitem{Lewy} {\sc H. Lewy,}
On the non-vanishing of the Jacobian in certain one-to-one mappings,
\textit{Bull. Amer. Math. Soc.,} {\bf 42}(1936), 689--692.


\bibitem{LX} {\sc X.~Y.~Liu,}
Bloch functions of several complex variables,
\textit{Pacific J. Math.}, {\bf 152}(1992), 347--363.

\bibitem{M3} {\sc M.~Mateljevi\'c,}
Quasiconformal and quasiregular harmonic analogues of Koebe's
theorem and applications,
\textit{Ann. Acad. Sci. Fenn. Math.}, {\bf 32}(2007),  301--315.

\bibitem{P} {\sc Ch.~Pommerenke,}
{\it Boundary behaviour of conformal maps,} Springer-Verlag Berlin
Heidelberg New York, 1992.


\bibitem{R}{\sc W.~Rudin,}
{\it Function theory in the unit ball of $\mathbb{C}^{n}$,}
Spring-Verlag, New York, Heidelberg, Berlin, 1980.

\bibitem{T} {\sc S.~Takahashi,}
Univalent mappings in several complex variables,
\textit{Ann. of Math.}, {\bf 53}(1951), 464--471.

\bibitem{W} {\sc H.~Wu,}
Normal families of holomorphic mappings,
\textit{Acta Math.}, {\bf 119}(1967), 193--233.

\end{thebibliography}
\end{document}